%\magnification=\magstep{1}
%\pagewidth{5.65in}
%\pageheight{7.5in}
\hoffset0.4 true in
\voffset-30pt
%\TagsOnRight
\input amstex
\documentstyle{amsppt}
\nologo
\topmatter
\title {An F. and M. Riesz theorem for planar vector fields}
\endtitle
\author {S. Berhanu and J. Hounie}
\endauthor
\address 
Department of Mathematics, Temple University, 
Philadelphia, PA 19122-6094, USA 
\endaddress
\email 
berhanu\@euclid.math.temple.edu
\endemail
\address 
Departamento de Matem\'atica, UFSCar,
13.565-905,  S\~ao Carlos, SP,  BRASIL
\endaddress
\email 
hounie\@ufscar.dm.br
\endemail
\thanks
 Work supported in part by CNPq, FINEP and FAPESP.
\endthanks
\def\ce{\Bbb C}
\def\erre{\Bbb R}
\def\ccinf{C^\infty_{c}}

\define\<{\langle}

\define\>{\rangle}

\redefine \|{\Vert}
\def\re{\hbox{\rm Re}\, }
\def\im{\hbox{\rm Im}\, }

\abstract
We prove that solutions of the homogeneous equation $Lu=0$, where $L$ is a locally integrable vector field with
smooth coefficients in two variables possess the F. and M. Riesz property.
That is,  if $\Omega$ is an open subset of the plane with smooth boundary,
$u\in C^1(\Omega)$ satisfies $Lu=0$ on $\Omega$, has  tempered growth at the
boundary, and its weak boundary value is a measure $\mu$, then $\mu$ is
absolutely continuous with respect to Lebesgue measure on the
noncharacteristic portion of  $\partial\Omega$.  \endabstract

 \endtopmatter
 \NoBlackBoxes
 \document
\smallskip
\heading
Introduction
\endheading
Consider a Borel measure $\mu$ defined on the boundary $\Bbb T$ of the unit circle $\Delta$ of the complex plane. A classical theorem proved in 1916 by  F.~and M.~Riesz states that if the Fourier coefficients of $\mu$ vanish for all negative integral values, i.e.,
$$
\widehat\mu(k)=\int_0^{2\pi}\exp(-2\pi i k\theta)\,d\mu(\theta)=0,\qquad k=-1,-2,\dots,\tag a
$$
then $\mu$ is absolutely continuous with respect to the Lebesgue measure $d\theta$. Condition (a) is equivalent to the existence of a holomorphic function $f(z)$ defined on $\Delta$ whose weak boundary value is $\mu$.

The F.~and M.~Riesz theorem has undergone  an extensive generalization  in
the last decades, mainly in two different directions: i) generalized
analytic function algebras, which has as a starting point the fact that
(a) means that $\mu$ is orthogonal to the algebra of continuous
functions on    $\Bbb T$ that extend holomorphically to $\Delta$; ii)
ordered groups, which emphasizes instead  the role of the group
structure of $\Bbb T$ in the classical result. Thus, although absolute
continuity with respect to Lebesgue measure is a local property (i.e.,
if each point has a neighborhood where it holds then it holds
everywhere), both directions focus on global objects. A remarkable exception is the paper  [B] in which the author uses  microlocal analysis to prove some generalizations of the theorem of F.~and M.~Riesz.   Among other things, in [B] it is shown that if a CR measure on a hypersurface of $\ce^n$ is the boundary value of a holomorphic function defined on a side, then it is absolutely continuous with respect to Lebesgue measure.

In view of Riemann's mapping theorem and the local character of the
conclusion, another way of stating the F.~and M.~Riesz theorem is to
say that if a holomorphic function $f(z)$ defined on a smoothly bounded
domain $D$ of the complex plane has tempered growth at the boundary and
its weak boundary value is a measure, then the measure is absolutely
continuous with respect to Lebesgue measure. If we regard holomorphic
functions as solutions of the homogeneous equation $\overline\partial
f=0$, it is natural to ask for which complex vector fields $L$ it is
possible to draw the same conclusion for solutions of the equation
$Lf=0$. In this paper we extend the F. and M. Riesz theorem to all
locally integrable, smooth complex vector fields in the plane for
smooth domains at the noncharacteristic part of the boundary. We
recall that a nowhere vanishing smooth vector field $$
L=a(x,y)\frac{\partial}{\partial x}+b(x,y)\frac{\partial}{\partial y}
$$ is said to be locally integrable in an open set $\Omega$ if each
$p\in \Omega$ is contained in a neighborhood which admits a smooth
function $Z$ with the properties that $LZ=0$ and the differential
$dZ\neq 0$. Examples of locally integrable vector fields include real
analytic vector fields and smooth, locally solvable vector fields. We
note however that the class of locally integrable vector fields is much
larger and refer to [T] on this subject.

 In his work [B], the author gives a microlocal criterion for the
absolute continuity of a measure analogous to (a) based on Uchiyama's
deep characterization of $BMO(\erre^n)$ [U]. Similarly, one of the main
steps in our generalization of the F.~and M.~Riesz theorem is Theorem
3.1 in section 3 which concerns the location of the wave front set of
the trace of a $C^1$ solution of a locally integrable vector field in
${\Bbb R}^n$. Here an important tool is the use of the FBI transform in
the fashion developed in [BCT] and [T]. However, since we apply this
result for $n=2$, in which case the trace lives in a one dimensional
boundary, we do not need to rely on Uchiyama's theorem and the
classical criterion (a) suffices. On the other hand, while in the
classical case and the generalizations in [B] the location of the wave
front set of the measure under consideration always satisfies a
restrictive hypothesis which leads to absolute continuity, this
restriction is not fulfilled in general by the trace of a solution of
an arbitrary locally integrable vector field even if the solution is
smooth (an example concerning a
vector field with real analytic coefficients is shown in Example 4.3).
Thus, we need to deal as well with points where the wave front set of
the measure may contain all directions; at those points the vector
field $L$ exhibits a behavior close to that of a real vector field (in
a sense made precise in Lemma 3.3) and absolute continuity may be
proved directly.

The paper is organized as follows. In section 1 we state our main
result (Theorem 1.1) and prove a lemma on the existence of  traces
on a noncharacteristic boundary for continuous solutions of an arbitrary
smooth complex vector field. In Section 2
we slightly extend one of the basic results in [B] with a new method of
proof  based on the FBI transform. In Section 3 we focus on locally integrable
vector fields and give a refined description of the location of the
wave front set of a solution's trace; this is a key step in the proof of
Theorem 1.1 which is also proved in this section. Finally, in section 4 we
present some examples.

\heading 1.
Statement of the main result
\endheading

The following theorem is the main result of this article. The existence of
the trace $bf=f(x,0)$ will be proved in Lemma 1.2 in this section.
 \proclaim{Theorem 1.1} Suppose $L=\frac{\partial}{\partial
t}+a(x,t)\frac{\partial}{\partial x} $ is smooth in a neighborhood $U$ of the
origin in the plane. Let $U_+=U\cap{\Bbb R}_+^{2}$, and suppose $f \in
C^1(U_+)$ satisfies $Lf=0$ in $U_+$ and for some integer $N$, $$
|f(x,t)|=O(t^{-N})\qquad\text{ as } t\to 0.
$$
Assume that $L$ is locally integrable in $U$. If the trace $bf=f(x,0)$ is a
measure, then it is absolutely continuous with respect to Lebesgue measure.
\endproclaim

 We begin with a lemma on the existence of a trace that improves Theorem 3.4 in
[B]. In the lemma, the vector field $L$ will not be assumed to be
locally integrable.
\proclaim{Lemma 1.2} Let $X\subseteq {\Bbb
R}^n $ be open, $U$ an open neighborhood  of $X\times\{0\}$ in
${\Bbb R}^{n+1}$, $U_+=U\cap{\Bbb R}_+^{n+1}$.  Let
$L=\frac{\partial}{\partial
t}+\sum_{j=1}^{n}a_j(x,t)\frac{\partial}{\partial x_j}$,
$a(x,t)\in C^{\infty}$ on $X\cup U_+$. Let $f$ be a continuous
function on $U_+$ such that $Lf\in L^{\infty}(U_+)$ and for some
$N\in {\Bbb N}$, $$
   |f(x,t)|=O(t^{-N}) \qquad \text{ as } t\to 0,
 $$
uniformly on compact subsets of $X$. Then $\lim_{t\to 0}f(x,t)=bf$
exists in ${\Bbb D}'(X)$. Furthermore, if $X\times (0,T]\subseteq
U_+$, then the distributions $\{f(.,t):0\leq t \leq T\}$ are
uniformly bounded in ${\Bbb D}'(X)$.
\endproclaim

\proclaim{Remark}
{\rm In [B] this Lemma is proved under the
additional assumption that $f$ is $C^1$ and that 
$$ 
|\partial_x f(x,t)|=O(t^{-N})\qquad \text{ as } t\to 0.
$$}
\endproclaim
\demo{Proof of Lemma 1.2} We will proceed as in [B] with some
modifications. Let $\phi \in C_0^{\infty}(X)$, and $T>0$ such that
$$ \text{ supp }\phi\times [0,T]\subseteq X\cup U_+ $$ For
$\epsilon \geq 0$ sufficiently small, set
$$L^{\epsilon}=\frac{\partial}{\partial
t}+\sum_{j=1}^{n}a_j(x,t+\epsilon)\frac{\partial}{\partial x_j}$$
Let $k\in {\Bbb N}$. We will choose
$\phi_0^{\epsilon},\dots,\phi_k^{\epsilon} \in
C^{\infty}(\overline{U_+})$ such that if $$
\Phi^{k,\epsilon}(x,t)=\sum_{j=0}^k\phi_j^{\epsilon}(x,t)\frac{t^j}{j!},
$$ then $$ (1)\quad \Phi^{k,\epsilon}(x,0)=\phi(x),\quad \text{
and }\quad (2)\quad |(L^{\epsilon})^{\ast}\Phi^{k,\epsilon}(x,t)|
\leq Ct^k,$$ where $C$ depends only on the derivatives of $\phi$
upto order $k+1$. In particular, $C$ will be independent of
$\epsilon$. Define $\phi_0^{\epsilon}(x,t)=\phi(x)$. For $j\geq
1$, write $$ L^{\epsilon}=\frac{\partial}{\partial
t}+Q^{\epsilon}(x,t,\frac{\partial}{\partial x}) ,$$ and define $$
\phi_j^{\epsilon}(x,t)=-\frac{\partial}{\partial t}
\phi_{j-1}^{\epsilon}(x,t)+(Q^{\epsilon})^{\ast}\phi_{j-1}^{\epsilon}$$
One easily checks that (1) and (2) above hold with these choices
of the $\phi_j^{\epsilon}$. We will next use the integration by
parts formula of the form $$\int u(x,T)w(x,T)dx-\int
u(x,0)w(x,0)dx=\int_0^T\int_{{\Bbb R}^n}(wPu-uP^{\ast}w)dxdt$$
which is valid for $P$ a vector field, $u$ and $w$ in $C^1({\Bbb
R}^n\times [0,T])$ and the $x-$support of $w$ contained in a
compact set in ${\Bbb R}^n$. Note that the $x$-support of $
\Phi^{k,\epsilon}(x,t)$ is contained in the support of $\phi(x)$.
Let $\psi \in C_0^{\infty}(B_1(0))$, where $B_1(0)$ denotes the
ball of radius $1$ centered at the origin in ${\Bbb R}^{n+1}$.
Assume $\int \int \psi dxdt=1$, and for $\delta >0$, let
$\psi_{\delta}(x,t)=\frac{1}{\delta
^{n+1}}\psi(\frac{x}{\delta},\frac{t}{\delta})$. For $\epsilon
>0$, set $f_{\epsilon}(x,t)=f(x,t+\epsilon)$. Observe that if
$\delta < \epsilon$, then the convolution
$f_{\epsilon}*\psi_{\delta}(x,t)$ is $C^{\infty}$ in the region
$t>0$. In the integration by parts formula above set
$u(x,t)=f_{\epsilon}*\psi_{\delta}(x,t)$, $w(x,t)=
\Phi^{k,\epsilon}(x,t)$ and $P=L^{\epsilon}$. We get: $$\align
  \int_Xf_{\epsilon}*\psi_{\delta}(x,0) \phi(x)dx &=\int_X f_{\epsilon}*\psi_{\delta}(x,T) \Phi^{k,\epsilon}(x,T)dx \\
                               &-\int_0^T\int_XL^{\epsilon}\left (f_{\epsilon}*\psi_{\delta}\right )\Phi^{k,\epsilon}dxdt \tag 1.1 \\
                               &+\int_0^T\int_X f_{\epsilon}*\psi_{\delta}(L^{\epsilon})^{\ast}\Phi^{k,\epsilon}dxdt
  \endalign $$
Fix $\epsilon >0$. Let $\delta \rightarrow 0^+$. Note that
$f_{\epsilon}*\psi_{\delta}(x,t)$ converges uniformly to
$f_{\epsilon}(x,t)$ on a neighborhood $W$ of $\text{ supp } \phi
\times [0,T]$. Hence in ${\Bbb D}'(W)$, $$L^{\epsilon}\left
(f_{\epsilon}*\psi_{\delta}\right )\rightarrow
L^{\epsilon}f_{\epsilon}$$ as $\delta \rightarrow 0^+$.
Moreover, $L^{\epsilon}f_{\epsilon}(x,t)=Lf(x,t+\epsilon)\in L^{\infty}$. Hence by
Friederichs' Lemma,$$L^{\epsilon}\left
(f_{\epsilon}*\psi_{\delta}\right )\rightarrow
L^{\epsilon}f_{\epsilon}$$ in $L^2(W)$ as $\delta \rightarrow
0^+$. We thus get $$\align
  \int_Xf(x,\epsilon) \phi(x)dx &=\int_X f(x,T+\epsilon) \Phi^{k,\epsilon}(x,T)dx \\
                               &-\int_0^T\int_XL^{\epsilon}f_{\epsilon}(x,t)\Phi^{k,\epsilon}(x,t)dxdt \tag 1.2 \\
                               &+\int_0^T\int_X f_{\epsilon}(x,t)(L^{\epsilon})^{\ast}\Phi^{k,\epsilon}(x,t)dxdt
  \endalign $$

  In the third integral on the right, we have
 $$|f_{\epsilon}(x,t)(L^{\epsilon})^{\ast}\Phi^{k,\epsilon}(x,t)|\leq Ct^{k-N} ,$$
 where $C$ depends only on the derivatives of $\phi$ upto order $k+1$.
 Choose $k=N+1$. By the dominated convergence theorem, as $\epsilon \to 0$,
 this third integral converges to
$\int_0^T\int_X fL^{\ast}\Phi^kdxdt$. In the second integral on the right,
note that since $Lf\in L^2(X\times (0,T)) $, as $\epsilon\to 0$, the
translates $L^{\epsilon}f_{\epsilon}=(Lf)_{\epsilon}\to Lf$ in $L^2$. We thus get  $$\langle
bf,\phi \rangle = \int_X f(x,T)\Phi^{k}(x,T)ds-\int_0^T\int_XLf\Phi^kdx dt+
\int_0^T\int_X fL^{\ast}\Phi^kdxdt ,$$ where $\Phi^k=\Phi^{k,0}$. From formula
(1.2), we also see that there is $C>0$ independent of $\epsilon$ such that $$
|\langle f(.,\epsilon),\phi \rangle |\leq C\sum_{|\alpha |\leq k+1}\Vert
\partial^{\alpha}\phi \Vert_{L^{\infty}}\tag 1.3$$ \enddemo

\heading
2. The FBI approach
\endheading

We will next present another proof of Theorem 3.5 in [B]. Our method
of proof is based on a variant of the  FBI transform  developed in [BCT] and [T]; this approach will also be used in the proof of Theorem 1.1 which concerns a
class of vector fields not covered in [B]. In our version, thanks to
Lemma 1.2, we will not assume that $\partial_xf(x,t)$ has a tempered
growth as $t\to 0+$. Note also that in Theorem 2.1 local integrability
of $L$ is not assumed.

\proclaim{Theorem 2.1} Let $X$, $U$, $U_+$ and $L=\frac{\partial}{\partial
t}+\sum_{j=1}^{n}a_j(x,t)\frac{\partial}{\partial x_j}$ be as in Lemma
1.2. Suppose $f\in C^1(U_+)$ satisfies $|f(x,t)|=O(t^{-N})$ for some
$N$ and $$|Lf(x,t)|=O(t^k),\quad k=1,2,...$$ uniformly on compact
subsets of $X$. Assume $$\partial_t^j a(x,0)=0\quad \forall j<l,\quad
\forall x\in X$$ and that $$\langle\partial_t^l
\im a(x_0,0),\xi^0\rangle >0\quad \text{ for some } x_0\in X,\quad
\xi^0\in{\Bbb R}^n.$$ Then $(x_0,\xi^0)\notin WF(bf)$
\endproclaim

\demo{Proof} Without loss of generality, we may assume that
$x_0=0$. Let $Z_1,\dots,Z_n$ be a complete set of smooth
approximate first integrals of $L$ near the origin in $U$
(see [T] for the existence of such). That is, $$LZ_j(x,t)=O(t^k),\quad
k=1,2,...\quad \text{ and } \qquad Z_j(x,0)=x_j,\quad 1\leq j\leq
n .$$ For $j=1,\dots,n$ let $M_j=\sum_{k=1}^n
b_{jk}(x,t)\frac{\partial}{\partial x_k}$ be vector fields
satisfying $$ M_jZ_k=\delta_j^k,\qquad [M_j,M_k]=0.$$ Note that
for each $j$, $$[M_j,L]=\sum_{s=1}^n c_{js}M_s \tag  2.1$$  where
each $c_{js}=O(t^k),$ $k=1,2,...$ Indeed, the latter can be seen
by expressing $[M_j,L]$ in terms of the basis
$\{L,M_1,\dots,M_n\}$ and applying both sides to the $n+1$
functions $\{t,Z_1,\dots,Z_n\}$. For any $C^1$ function $g$,
observe that the differential $$dg=\sum_{k=1}^n
M_k(g)dZ_k+(Lg-\sum_{k=1}^n M_k(g)LZ_k)dt \tag  2.2$$ This is
verified by evaluating each side at the basis vector fields
$\{L,M_1,\dots,M_n\}$. Using (2.2) we get:
$$d(gdZ_1\wedge\dots\wedge
dZ_n)=(Lg-\sum_{k=1}^nM_k(g)LZ_k)dt\wedge dZ_1\wedge\dots\wedge
dZ_n \tag  2.3$$ For $\xi \in {\Bbb R}^n,\quad s\in {\Bbb R}^n$,
let $$E(s,\xi,x,t)=i\xi\cdot (s-Z(x,t))-|\xi|(s-Z(x,t))^2 ,$$
where for $w\in {\Bbb C}^n$, we write $w^2=\sum_{j=1}^n w_j^2$.
Let $B$ denote a small ball centered at $0$ in ${\Bbb R}^n$ and
$\phi \in C_0^{\infty}(B)$, $\phi \equiv 1$ near  the origin. We
will apply (2.3) to the function $$
g(s,\xi,x,t)=\phi(x)f(x,t)e^{E(s,\xi,x,t)} $$ where $(s,\xi)$ are
parameters. We get: $$ d(gdZ)=\{L(\phi f)+(\phi
f)LE-\sum_{k=1}^n(M_k(\phi f)+\phi f(M_kE))LZ_k\}e^Edt\wedge
dZ,\tag  2.4$$ where $dZ=dZ_1\wedge \dots \wedge dZ_n$. Next by
Stokes theorem we have, for $t_1>0$ small: $$\int_B
g(s,\xi,x,0)dx=\int_B g(s,\xi,x,t_1)d_xZ(x,t_1)+\int_0^{t_1}\int_B
d(gdZ) \tag  2.5$$ We will estimate the two integrals on the right
in (2.5). Write
   $$Z=(Z_1,\dots,Z_n)=x+t\Psi(x,t), \text{ and } \Psi=\Psi_1 +i\Psi_2.$$
   Since the $Z_j$ are approximate solutions of $L$, we have
$$\Psi+t\Psi_t+(I+t\Psi_x)\cdot a=O(t^k),\quad k=1,2...$$
and hence
$$\partial_t^j\Psi(x,0)=0 ,\quad j<l\text{ and } \langle
\partial_t^l\Psi_2(x,0),\xi^0\rangle <0 \tag  2.6$$ for $x$ in a neighborhood
$V$ of $\overline B$ (after shrinking $B$, if necessary).
Observe that
$$
\re E(s,\xi,x,t)=t\xi\cdot\Psi_2(x,t)-|\xi|((s-x-t\Psi_1)^2
-t^2\Psi_2(x,t)^2)
$$
Because of (2.6), continuity and homogeneity in $\xi$, we can get
$c_1>0$ such that $$\re E(s,\xi,x,t)\leq -c_1|\xi|t^{l+1},\quad \text{
for } x\in V,\quad 0\leq t\leq t_1 $$ $$s\in {\Bbb R}^n \text{ and }
\xi \text { in a conic nbhd } \Gamma \text{ of } \xi^0. \tag 2.7$$
Going back to the integrals in (2.5), we clearly have

$$
\left|\int_Bg(s,\xi,x,t_1)d_xZ(x,t_1)\right|\leq e^{-c_2|\xi|},
$$
for some $c_2>0$, for $s\in {\Bbb R}^n$ and $\xi \in \Gamma$.  To
estimate $\int_0^{t_1}\int_B d(gdZ)$, we use (2.4) and look at each
term that appears there. We first consider the term $L(\phi f)e^E$. For
any $k$,

$$
|\phi(Lf)e^E|\leq C_kt^{lk}e^{-c_1t^l|\xi|}\leq \frac{C_k'}{|\xi|^k}
$$
Moreover, the $x$-integral
$$
\int_B (L\phi)fe^EdZ=\langle f(.,t),(L\phi)e^E\rangle
$$
can be estimated using Lemma 1.2.  Accordingly, after decreasing $t_1$,
we can get $\delta >0$ such that if $|s|\leq \delta$ and $\xi\in
\Gamma$,

$$
|\langle
f(.,t),(L\phi)e^E\rangle | \leq C\left|\sum_{|\alpha| \leq N+1}
\partial_x^{\alpha}((L\phi) e^E)\right|_{L^{\infty}}\leq Ce^{-c|\xi|}
$$
for some constants $c, C>0$. In the latter, we have used the constancy
of $\phi$ near $0$. It follows that the integral
$$\int_B\int_0^{t_1}L(\phi f)e^Edt\wedge dZ$$ decays rapidly in $\xi$.
The term $(\phi f)LEe^E$ is estimated using the fact that for any $k$,
$|LE|\leq c_kt^k|\xi|$ for some constant $c_k$ and that $|e^E|\leq
e^{-c_1t^l|\xi|}$. This shows that
$$
\int_B\int_0^{t_1}(\phi f)LEe^E dt\wedge dZ
$$
decays rapidly in $\xi$. The integral of $\phi f(M_kE)LZ_ke^E$ is
estimated likewise. To estimate the integral of $(M_k(\phi f))LZ_ke^E$,
we first integrate in $x$ and apply Lemma 1.2 again. Indeed, the Lemma
also applies to the weak derivative $M_k(\phi f)$. Thus

$$
\int_B\int_0^{t_1}d(gdZ)
$$
has a rapid decay in $\xi$, and going back to (2.5), we have shown:
$$
F(s,\xi) = \int_Be^{i\xi
\cdot(s-x)-|\xi|(s-x)^2}\phi(x) f(x,0)dx \tag 2.8
$$
decays rapidly for $|s|\leq \delta$ in ${\Bbb R}^n$ and $\xi$ in a conic neighborhood
$\Gamma $ of $\xi^0$. The function $F(s,\xi)$ is the FBI transform
(see [BCT]) of the distribution $\phi (x)f(x,0)$. To conclude the
proof, we will exploit the inversion formula for the FBI, namely, $$
\phi(x)f(x,0)=\lim_{\epsilon\to 0^+}c_n\iint e^{i(x-s)\cdot \xi
-\epsilon |\xi|^2}F(s,\xi)|\xi|^{\frac{n}{2}}dsd\xi \tag 2.9$$ where
$c_n$ is a dimensional constant. Assume now that $\phi(x)$ is supported
in the ball centered at the origin with radius $M$. We will study the
inversion integral in (2.9) by writing it as a sum of three pieces :
$I_1(\epsilon)$, $I_2(\epsilon)$, and $I_3(\epsilon)$. The first piece
consists of integration over the region $\{(\xi,s):|s|\geq 2M\}$. In
the second piece we integrate over $\{(\xi,s):\delta \leq |s| < 2M\}$,
and in the third piece over $\{(\xi,s):|s|\leq\delta\}$. For the
integral $I_1(\epsilon)$, after integrating in $s$, one gets an
exponential decay in $\xi$ independent of $\epsilon$, and hence
$\lim_{\epsilon\to 0^+} I_1(\epsilon)$ is in fact a holomorphic
function near the origin in ${\Bbb C}^n$. To study the second piece, we
write it as $$I_2(\epsilon)=c_n\int_{\{(y,\xi,s):\delta\leq
|s|<2M\}}e^{i(x-y)\cdot\xi-|\xi|(s-y)^2-\epsilon
|\xi|^2}\phi(y)f(y,0)|\xi|^{\frac{n}{2}}dydsd\xi$$ We will use the
holomorphic function $\langle \zeta \rangle=(\zeta_1^2+\cdots
+\zeta_n^2)^{\frac{1}{2}}$ where we take the principal branch of the
square root in the region $|\im \zeta |<|\re \zeta|$. Observe that this
function is a holomorphic extension of $|\xi|$ away from the origin. In
the $\xi$ integration above, we can deform the contour to the image of
$$\zeta(\xi)=\xi+i\beta(x-y)|\xi|$$ where $\beta$ is chosen
sufficiently small. In particular, we choose $\beta$ so that when $x$
varies near the origin and $y$ stays in the support of $\phi$, then
$|\im \zeta(\xi)| < |\re \zeta(\xi)|$, away from $\xi=0$. In the
integrand of $I_2(\epsilon)$, if $|x|\leq \frac{\delta}{4}$, we get an
exponential decay independent of $\epsilon$. It follows that this piece
is also holomorphic near the origin in ${\Bbb C}^n$ after setting
$\epsilon =0$. Finally, for the third piece, let
$\Gamma_1,\dots,\Gamma_m$ be convex cones such that with
$\Gamma_0=\Gamma$, $${\Bbb R}^n=\bigcup_{j=0}^{m}\Gamma_j,$$ and for
each $j\geq 1$ there exists a vector $v_j$ satisfying $v_j\cdot
\overline{\Gamma_j}>0$ and $v_j\cdot \xi^0<0$. We now write
$$I_3(\epsilon)=\sum_{j=0}^m K_j(\epsilon),$$ where $K_j$ equals the
integral over $\Gamma _j$. The decay in the FBI established in (2.8)
shows us that $K_0$ is a smooth function even after setting
$\epsilon=0$. Each of the remaining functions $K_j$, after setting
$\epsilon =0$, is a boundary value of a tempered holomorphic function
in a wedge whose inner product with $\xi^0$ is negative. Hence
$$(0,\xi^0) \notin WF_a (K_j(0+)),$$ where $WF_a$ denotes the analytic
wave front set (see [S]). The latter implies that
$$
(0,\xi^0) \notin WF (K_j(0+)).
$$
Indeed, as is well known, a distribution $u$ is
microlocally analytic (resp. smooth) at a covector $\gamma$ iff there is an
analytic (resp. smooth) pseudodifferential operator $P$ elliptic at
$\gamma$ such that $Pu$ is analytic (resp. smooth). We have thus proved that
$(0,\xi^0)\notin WF (f(x,0))$.

\enddemo

\heading 3.
Proof of Theorem 1.1 and auxiliary results
\endheading

We next wish to get a better description of the wave front set of the trace
of  a solution when the vector field in question is locally
integrable.
\newline 
We consider a smooth vector field $ L=X+iY$
where $X$ and $Y$ are real vector fields defined in a neighborhood
$U$ of the origin. Let $\Sigma$ be an embedded hypersurface
through the origin in $U$ dividing the set $U$ into two regions,
$U^+$ and $U^-$ where $U ^+$ denotes the region  towards which $X$
is pointing. We will consider a function $f\in C^1(U^+)$ that
satisfies $Lf=0$ on $U^+$ and grows in a tempered fashion as
$p\mapsto \Sigma$. We assume that $L$ is noncharacteristic on
$\Sigma$ which means (after multiplying $L$ by $i$ if necessary)
that $X$ is noncharacteristic. Our considerations will be local
and so after an appropriate choice of local coordinates $(x,t)$
and multiplication of $L$ by a nonvanishing factor, the vector
field is given by $$
 L=\frac{\partial}{\partial t}+\sum_{j=1}^{n}a_j(x,t)\frac{\partial}{\partial
x_j} \tag  3.1
$$
and $\Sigma$ and $U^+$ are given by $t=0$ and $t>0$
respectively.

We will need to consider the integral curve $(-\epsilon,\epsilon) \owns
s \mapsto \gamma(s)$ of $X$ that passes though the origin, i.e.,
$\gamma'(s)= X\circ\gamma(s)$, $\gamma(0)=0$. It is clear that for
small $\epsilon>0$ and $|s|<\epsilon$, $\gamma(s)\in U^+$ if and only
if $s>0$, so $\gamma((-\epsilon,\epsilon))\cap U^+=
\gamma((0,\epsilon))$. To simplify the notation we will simply write
$\gamma^+$ to denote $\gamma((0,\epsilon))$.

\proclaim{Theorem 3.1}
Let  $L=\frac{\partial}{\partial
t}+\sum_{j=1}^{n}a_j(x,t)\frac{\partial}{\partial x_j}$ be locally integrable.
\roster
\item"i)"
Suppose $f\in C^1(U_+)$ satisfies $|f(x,t)|=O(t^{-N})$ for some $N$ and  $$
Lf(x,t)=0, \quad (x,t)\in U^+.
$$
Assume that there is a sequence $p_k\in \gamma^+$, $p_k\to 0$ such that for
each $k=1,2,\dots$,  $X(p_k)$ and $Y(p_k)$ are linearly
independent. Then there exists a unit vector $v$ such that
$$
\xi^0\in{\Bbb R}^n,\quad v\cdot\xi^0>0\quad\Longrightarrow
(0,\xi^0)\notin WF(bf).
$$
In particular, the wave front set of $bf$ at the origin is contained in a closed half-space.
\item"ii)"
Conversely, if $X(p)$ and $Y(p)$ are linearly dependent for all $p\in\gamma^+$ there exists a neighborhood $V\subset U$ of the origin and a function $f\in C^1(V_+)\cap C^0(\overline{V_+})$ such that $Lf=0$ and $(0,\xi)\in WF(bf)$ for all $\xi\in\erre^n\setminus\{0\}$.
\endroster

\endproclaim
\demo{Proof} i) Let $Z_1,\dots,Z_n$ be a complete set of smooth
 first integrals of $L$ near the origin in $U$ (see [T] on the subject of locally integrable structures). That is,
$$ LZ_j(x,t)=0, \quad k=1,\dots,n,\quad dZ_1\wedge\cdots\wedge
dZ_n(0,0)\not=0, $$ and  choose new  local coordinates $(x,t)$ in
which the $Z_j$'s may be written as $$ Z_j(x,t)=x_j+i\Phi_j(x,t),
\quad k=1,\dots,n, $$ with $\Phi(0,0)=0$ , $\Phi_x(0,0)=0$ and
$\Phi_{xx}(0,0)=0$ ($\Phi=(\Phi_1,\dots,\Phi_n)$).

For $j=1,\dots,n$ let $M_j=\sum_{k=1}^n
b_{jk}(x,t)\frac{\partial}{\partial x_k}$ be vector fields satisfying

$$
M_jZ_k=\delta_j^k,\qquad [M_j,M_k]=0.
$$
It is readily checked that for each $j=1,\dots,n$,
$$
[M_j,L]=0.\tag 3.2
$$
 For any $C^1$ function $g$,  the differential may be expressed as
$$ dg=Lg\,dt +\sum_{k=1}^n M_kg\,dZ_k\tag 3.3$$ Using (3.3) we
get: $$d(gdZ_1\wedge\dots\wedge dZ_n)=Lg\,dt\wedge
dZ_1\wedge\dots\wedge dZ_n \tag 3.4$$
 For $\zeta\in {\Bbb C}^n,\quad z\in {\Bbb C}^n$, let $$
E(z,\zeta,x,t)=i\zeta\cdot (z-Z(x,t))-\kappa\langle \zeta
\rangle(z-Z(x,t))^2 , $$ Let $B$ denote a small ball centered at $0$ of
radius $r$ in ${\Bbb R}^n$ and $\phi \in C_0^{\infty}(B)$, $\phi \equiv
1$ for $|x|\le r/2$, the precise value of $r$ as well as the value of
the positive parameter $\kappa$ in the definition of $E$ will be
determined later. We will apply (3.4) to the function

$$
g(z,\zeta,x,t)=\phi(x)f(x,t)e^{E(z,\zeta,x,t)} $$ where
$(z,\zeta)$ are parameters. We get:
$$
d(gdZ)=f\,L\phi\,e^Edt\wedge dZ,\tag 3.5
$$
where $dZ=dZ_1\wedge
\dots \wedge dZ_n$. Next by Stokes theorem we have, for $t_1>0$
small: $$\int_B g(z,\zeta,x,0)d_xZ(x,0)=\int_B
g(z,\zeta,x,t_1)d_xZ(x,t_1)+\int_0^{t_1}\int_B d(gdZ) \tag 3.6 $$
We will estimate the two integrals on the right in (3.6) and our
aim is to show that for $x$ and $z$ close to the origin  in real
and complex space respectively, both decay exponentially as $\zeta
\to\infty$ in a conic neighborhood of  $\xi^0$. Write $$
Z=(Z_1,\dots,Z_n)=x+i\Phi(x,t), \quad\Phi=(\Phi_1,\dots,\Phi_n).
$$ Observe that, assuming without loss of generality that
$|\xi^0|=1$, $$ \re
E(0,\xi^0,x,t)=\Phi(x,t)\cdot\xi^0-\kappa(|x|^2-|\Phi(x,t)|^2). $$
Our main task will be to determine convenient  values of $t_1$,
$\kappa$ and $r$ such that for some $\gamma>0$ \roster
\item $\re E(0,\xi^0,x,t_1)\le-\gamma$ for $|x|\le r$ ;

\item $\re E(0,\xi^0,x,t)\le -\gamma$ for $0\le t\le t_1$ and $r/2\le|x|\le r$.
\endroster
The assumptions on $\Phi$ allow us to write

$$ \Phi(x,t)=\Phi(0,t)+e(x,t),\qquad |e(x,t)|\leq A|xt|+B|x|^2
\tag  3.7$$ for some positive constants $A$ and $B$. We may
assume that $\Phi_t(0,0)=0$, otherwise the result we want to prove
would follow from Theorem  2.1. Hence, we may assume that the
quotient $|\Phi(0,t)|/t^2\le C$ for $(0,t)\in U^+$. In order to
find the vector $v$ mentioned in the statement of the theorem we
will need
\proclaim{Lemma 3.2} There exist a sequence $t_k\searrow
0$ such that \roster \item  $\Phi(0,t_k)\not=0$;

\item  $|\Phi(0,t)|\le |\Phi(0,t_k)|$ for $0\le t\le t_k$;

\item $\displaystyle\lim_{t_k\to0}\Phi(0,t_k)/|\Phi(0,t_k)|=-v$

\endroster
\endproclaim
We will postpone the proof of Lemma 3.2  and continue  our reasoning with $v$
given by (3). We have  $\Phi(0,t_k)+ |\Phi(0,t_k)| \,\,v=o(|\Phi(0,t_k)|)$. We
recall that by hypothesis $\xi^0\cdot v>0$. Hence,  $$ \aligned
\Phi(0,t_k)\cdot\xi^0
           &=-|\Phi(0,t_k)|\,\, v\cdot\xi^0+o(|\Phi(0,t_k)|)\\
           &<-|\Phi(0,t_k)| \,\,v\cdot\xi^0/2=-c|\Phi(0,t_k)|,\\
\endaligned
$$
for $t_k$ small and $0<c<1$. We now take $r=\alpha|\Phi(0,t_k)|/t_k$, with $\alpha$ and $t_k$ small to be chosen later.
Hence, for $|x|\le r$ and $0\le t\le t_k$, we can choose $\alpha$ small
enough (depending on $A$,  $B$ and $C$ but not on $t_k$)  so that $$ \aligned
|e(x,t)|&\le A\alpha\,|\Phi(0,t_k)|\,{t\over t_k}+  B\alpha^2\,{|\Phi(0,t_k)|\over t_k^2}\,|\Phi(0,t_k)|\\
         &\le c{|\Phi(0,t_k)|\over 2}.  \\
\endaligned
\tag  3.8$$ This implies that on the support of $\phi(x)$ we have
$$ -(1+c)|\Phi(0,t_k)|\leq \Phi(x,t_k)\cdot\xi^0\leq - {c\over2}
|\Phi(0,t_k)|. $$ Let $\kappa=\epsilon/|\Phi(0,t_k)|$. A
consequence of (3.7), (3.8) and the fact that $|\Phi(0,t)|\le
|\Phi(0,t_k)|$ for $0\le t\le t_k$ is $$ \aligned |\Phi(x,t)|&\le
(1+c)|\Phi(0,t_k)|\\ |\Phi(x,t)|^2&\le (1+c)^2|\Phi(0,t_k)|^2\\
\kappa|\Phi(x,t)|^2&\le \epsilon(1+c)^2|\Phi(0,t_k)|\\
\endaligned
\tag  3.9
$$
 for $x$ in the support of $\phi(x)$ and $0\le t\le t_k$. Choosing
$\epsilon=c/(4(1+c)^2)$ (thus, independent of $t_k$), we get, on the
support of $\phi(x)$, $$
 \Phi(x,t_k)\cdot\xi^0+\kappa |\Phi(x,t_k)|^2\leq -{c\over2}|\Phi(0,t_k)|+\epsilon(1+c)^2|\Phi(0,t_k)|\le -{c\over4}
|\Phi(0,t_k)|$$
 which leads to an exponential decay in the first integral on the right
of (3.6) for $z$ complex near $0$ and $\zeta$ in a complex conic
neighborhood of $\xi^0$, as soon as we replace $t_1$ by $t_k$. For the
second integral, note that for $0\leq t\leq t_k$ and $x$ in the support
of $\phi$, we may invoke again (3.9) to estimate the size of
$|\Phi(x,t)|$ and $\kappa |\Phi(x,t)|^2$ which gives, in view of the
previous choice of $\epsilon$, 
$$
|\Phi(x,t)|+\kappa |\Phi(x,t)|^2\leq
(1+c)|\Phi(0,t_k)|+ {c\over4}|\Phi(0,t_k)|\le (1+2c)|\Phi(0,t_k)| $$
while on the support of $L\phi$, $|x|\ge r/2=\alpha |\Phi(0,t_k)|
/2t_k$ so $$ \kappa |x|^2\geq
{\epsilon\alpha^2|\Phi(0,t_k)|\over4t_k^2} $$ and $$
\Phi(x,t)\cdot\xi^0-\kappa(|x|^2-|\Phi(x,t)|^2)\le
    (1+2c-{\epsilon\alpha^2\over4t_k^2}) |\Phi(0,t_k)|.
$$ 
 Hence, if $t_k$ is chosen  sufficiently small, we also get exponential decay for the second integral on the right hand side of (3.6) with $t_1$ replaced by $t_k$.
\newline 
We have thus shown that the function
 $$
F(z,\zeta)=\int_Be^{E(z,\zeta,x,0)}\phi(x)f(x,0)d_xZ(x,0)$$
satisfies an exponential decay of the form $$ |F(z,\zeta)|\leq
Ce^{-R|\zeta|}$$ for $z$ near $0$ in ${\Bbb C}^n$ and $\zeta $ in
a complex conic neighborhood of $\xi^0$ in ${\Bbb C}^n$. In
particular, since $Z(0,0)=0$ and $d_xZ(0,0)$ is the identity
matrix, the function $$ G(x,\xi)=F(Z(x),(Z_x(x)^{-1})^t\xi) $$ has
an exponential decay for $(x,\xi)$ in a real conic neighborhood of
$(0,\xi^0)$. By Theorem 2.2 in [BC], it follows that
$(0,\xi^0)\notin WF(bf)$.

We now return to the proof of Lemma 3.2; it is here that  we use the fact
that $X$ and $Y$ are linearly independent on a sequence $p_k\in\gamma^+$ that
approaches the origin. We will show that $\Phi(0,t)$
cannot vanish identically on any interval $(0,\epsilon')$. Let us write
$L=\partial_t+a\cdot\partial_x$, $Z=x+i\Phi$, $Z_x=I+i\,^t\Phi_x$ and recall
that $^t\Phi_x$ has small norm for $(x,t)$ close to $0$. Now $LZ=0$ leads to
$a=-i(I+i\,^t\Phi_x)^{-1}  \, \Phi_t$. If $\Phi(0,t)$ vanishes identically on
$[0,\epsilon']$ we will have, for those values of $t$, that $\Phi_t(0,t)=0$,
$a(0,t)=0$, and $Y(0,t)=\im a(0,t)=0$. Furthermore, $X(0,t)=\partial_t$ for
$0<t< \epsilon'$, showing that $\gamma(s)=(0,\dots,0,s)$ for $0<s< \epsilon'$.
Thus, $X(\gamma(s))$ and $Y(\gamma(s))$ are linearly dependent for $0<s<
\epsilon'$, a contradiction. Therfore, there exists a sequence $s_k\searrow0$
such that  $|\Phi(0,s_k)|>0$ and since $\Phi(0,0)=0$ there is another sequence
 $t_k\searrow0$ satisfying (1) and (2),  which in turn possesses a subsequence
that satisfies (1), (2) and (3).

 ii) Consider as before a complete set of smooth  first integrals of $L=X+iY$ defined in neighborhood
 $V\subset U$ of the origin, $Z_1,\dots,Z_n$, $ LZ_j=0$, $k=1,\dots,n$,
$dZ_1\wedge \cdots\wedge dZ_n(0,0)\not=0$ , and  local coordinates $(x,t)$ in
which the $Z_j$'s have the form $ Z_j(x,t)=x_j+i\Phi_j(x,t)$, $k=1,\dots,n$,
with $\Phi(0,0)=0$, $\Phi_x(0,0)=0$ and $\Phi_{xx}(0,0)=0$. Since $Y$ is
proportional to $X$ along $\gamma^+$ it follows that $XZ_j=0$ on $\gamma^+$;
therefore $X\re Z_j=0$ on $\gamma^+$ which in the coordinates $(x,t)$ implies
that $x_j(\gamma^+)=0$, $j=1,\dots,n$. Since $Z_j$ vanishes on $\gamma^+$ we
conclude that   $\Phi_j(0,t)=0$ for $0\le t<t_0$, for some $t_0>0$ and any
$j=1,\dots,n$. This shows that $Z=(Z_1,\dots,Z_n)$ maps
$\{0\}\times[0,t_0)\subset\erre^n\times\erre$ into $\{0\}\subset\ce^n$.

Let us denote by  $F(\zeta)=\langle \zeta \rangle=(\zeta_1^2+\cdots
+\zeta_n^2)^{\frac{1}{2}}$  the holomorphic function $\langle \zeta \rangle=(\zeta_1^2+\cdots
+\zeta_n^2)^{\frac{1}{2}}$ where we take the principal branch of the square root in the region $|\im \zeta |<|\re \zeta|$ and set $F(0)=0$.
We also set $F_\epsilon(\zeta)=(\zeta_1^2+\cdots +\zeta_n^2+\epsilon)^{\frac{1}{2}}$, $\epsilon>0$.
Shrinking $V$ if necessary we assume that $Z(V_+)=V_++i\Phi(V_+)$ is contained in  $|\im \zeta |\le|\re \zeta|/2$ and the 
composition $u(x,t)=F(Z(x,t))$ is well defined and continuous. Approximating
$u(x,t)$ by the smooth functions $u_\epsilon(x,t)=F_\epsilon(Z(x,t))$,
$\epsilon\searrow0$ that satisfy  $Lu_\epsilon=0$ we see that $u$ satifies the
homogeneuous equation $Lu=0$ in the sense of distributions in $V_+$ and so
does $f(x,t)=u(x,t)^3$. A moment's reflection shows that $f\in C^1(V^+)\cap
C^0(\overline{V^+})$. Furthermore,
$Z_1^2(x,0)+\cdots+Z_n^2(x,0)=\lambda(x)\,|x|^2$ where $\lambda$ is a smooth
complex-valued function that does not vanish in a neighborhood of the origin,
so $f(x,0)=\lambda(x)^{3/2}|x|^3$ and the wave front sets of $f(x,0)$ and
$x\mapsto|x|^3$ coincide for small values of $x$.  The wave front set of
$|x|^3$ is precisely $\{0\}\times(\erre^n\setminus\{0\})$ because $|x|^3$ is
smooth except at the origin and it is invariant under rotations.

\enddemo

In the proof of Theorem 1.1 we will also need the following lemma on measures
which arise as traces of homogeneous solutions of vector fields.

\proclaim {Lemma 3.3} Let $$ L={\partial\over\partial
t}+i\sum_{j=1}^n b_j(x,t){\partial\over\partial x_j} $$ be smooth
on a neighborhood $U=B(0,a)\times(-T,T)$ of the origin in
$\erre^{n+1}$ with $B(0,a)=\{x\in\erre^n:\,\,\,|x|<a\}$. We will
assume that the coefficients $b_j(x,t)$, $j=1,\dots,n$ are real
and that all of them vanish on $F\times[0,T)$,  where $F\subset
B(0,a)$ is a closed set.   Assume that  $f\in C^1(U^+)$ satisfies
$Lf=0$ on $U^+=C^1(B(0,a)\times(0,T))$, has tempered growth as
$t\searrow0$ and  its boundary value $bf(x)=f(x,0)$ is a Radon
measure $\mu$. Then the restriction $\mu_F$  of $\mu$ to $F$
defined on Borel sets $X\subset B(0,a)$  by $\mu_F(X)=\mu(X\cap
F)$ is absolutely continuous with respect to Lebesgue measure.
\endproclaim
\demo{Proof}
 If $\tilde x$ is an arbitrary point in $F$ we may write
$$
b_j(x,t)=\sum_{k=1}^n (x_k-\tilde x_k)\beta_{jk}(x,\tilde x,t) \tag  3.10
$$
 with $\beta_{jk}(x,\tilde x,t)$ real and smooth.  Recall that for any
$\phi\in\ccinf(-a,a)$ we have $$ \aligned \< \mu,\phi \> &= \int
f(x,T)\Phi^{k}(x,T)ds+ \int_0^T\int_{-a}^a f(x,t)\,L^t\Phi^k(x,t)\,dxdt
\\ \Phi^k(x,t)&=\sum_{j=0}^k\phi_j (x,t)\frac{t^j}{j!}\\ \endaligned
\tag 3.11$$ where $\phi_0(x,t)=\phi(x)$, $$
\phi_j(x,t)=-\frac{\partial}{\partial t}
\phi_{j-1}^{\epsilon}(x,t)-\sum_{s,\ell=1}^n\frac{\partial}{\partial
x_s} (x_\ell-\tilde x_\ell)\beta_{j\ell}(x,\tilde x,t) \phi_{j-1}
(x,t), \quad j=1,\dots,k, $$ and $k$ is a convenient and fixed positive
integer. We may as well write $$ \Phi^k(x,t)=A(x,t,D_x)\phi(x)\tag 3.12
$$ where $A(x,t,D_x)=\sum_{|\alpha|\le k}\,a_\alpha(x,t)D^\alpha_x$ is
a linear differential operator of order $k$ in the $x$ variables with
coefficients depending smoothly on $t$.  The coefficients $a_\alpha$
are obtained from the coefficients $b_j(x,t)$ of $L$ by means of
algebraic operations and differentiations with respect to $x$ and $t$.
The key observation is that (3.10) implies that, given any point
$\tilde x\in F$, $A(x,t,D_x)$ may be written as $$
 A(x,t,D_x)=\sum_{|\alpha|\le k}\sum_{\ell=1}^n\,A_{\alpha\ell}(x,\tilde
x,t)\,((x_\ell-\tilde x_\ell) D_x)^\alpha. \tag 3.13$$
 Notice that $|A_{\alpha\ell}(x,\tilde x,t)|\le C$, for $x\in B(0,a)$,
$\tilde x\in F$, $t\in[0,T)$, $|\alpha|\le k$, and $\ell=1,\dots, n$
because the coefficients of $L$ have uniformly bounded derivatives on
$B(0,a)$. Hence, we obtain from (3.11), (3.12) and (3.13) the estimate
$$ \left | \int f(x,T)\Phi^k(x,T)dx\right | \le C\,\sum_{|\alpha|\le k+1}
\int_{B(0,a)}\,d(x,F)^{|\alpha|}\,|D_x^\alpha\phi(x)|\,dx,\tag 3.14$$
where $d(x,F)=\inf_{\tilde x\in F} |x-\tilde x|$. 
We next consider the second integral on the right in (3.11). We recall from
the proof of Lemma 1.2 that 
$$
|L^t\Phi ^k(x,t)|\leq Ct^k
$$
We need to examine this inequality more closely. We will first show that for
any $j$,
$$
L^t(\Phi ^j)=\frac{\phi_{j+1}}{j!}t^j \tag 3.15
$$
To see this, note that  (3.15) holds for $ j=0$ from the definition of
$\phi_1$. To proceed by induction, assume (3.15) for $j\leq m$. Then 
$$\align
 L^t(\Phi ^{m+1})&=L^t(\Phi ^m)+L^t\left
(\frac{\phi_{m+1}}{(m+1)!}t^{m+1}\right )\\
&=\frac{\phi_{m+1}}{m!}t^m+L^t\left (\frac{\phi_{m+1}}{(m+1)!}t^{m+1}\right )\\
&=\frac{L^t(\phi_{m+1})}{(m+1)!}t^{m+1}\\
&=\frac{\phi_{m+2}}{(m+1)!}t^{m+1}
\endalign
$$
This proves (3.15). Next we observe that since the coefficients $b_j(x,t)$
vanish on $F\times [0,T]$, each $\phi_j$ has the form
$$
\phi_j(x,t)=\sum_{|\alpha |\leq j}c_{\alpha}(x,t)D_x^{\alpha}\phi(x) 
$$
where the $c_{\alpha}$ are smooth and satisfy the estimate
$$
|c_{\alpha}|\leq Cd(x,F)^{|\alpha|} \tag 3.16
$$
The form (3.16) is clearly valid for $\phi_0=\phi$. Assume it is valid for
$\phi_j$. Then it will also be valid for $\phi_{j+1}$ since by definition,
$\phi_{j+1}=L^t\phi_j$. If we now choose $k=N+1$, (3.14), (3.15) and
(3.16) imply  that 
$$
\align
\left |\int_0^T\int_{-a}^af(x,t)L^t\Phi ^k(x,t)dxdt\right |&\leq
\int_0^T\int_{-a}^a|f(x,t)|\frac{\phi_{k+1}(x,t)}{k!}t^kdxdt\\
&\leq  C\int_0^T\int_{-a}^a|\phi_{k+1}(x,t)|dxdt\\
&\leq C\sum_{|\alpha|\leq
k+1}\int_{-a}^ad(x,F)^{|\alpha|}|D_x^{\alpha}\phi(x)|dx\tag 3.16
\endalign
$$
Thus the second integral on the right hand side of (3.11) also satisfies an
estimate of the kind in (3.14). Consider now a compact
subset $K\subset F$ with Lebesgue measure $|K|=0$ and choose a sequence
$0\le\phi_\epsilon(x)\le 1\,\in\ccinf(B(0,a))$, $\epsilon\to 0$, such
that i) $\phi_\epsilon(x)=1$ for all $x\in K$; ii) $\phi_\epsilon(x)=0$
if $d(x,K)>\epsilon$; iii) $|D^\alpha_x\phi_\epsilon(x)|\le
C_\alpha\epsilon^{-|\alpha|}$.  Note that $\phi_\epsilon(x)$ converges
pointwise to the characteristic function of $K$ as $\epsilon\to 0$
while $D^\alpha\phi_\epsilon(x)\to0$ pointwise if $|\alpha|>0$.  Let
$\psi\in\ccinf(B(0,a))$ and use (3.14) and (3.16) with $\phi=\phi_\epsilon\psi$
keeping in mind the trivial estimate $d(x,F)\le d(x,K)$. By the
dominated convergence theorem,
 $\<\mu,\phi_\epsilon\psi\>\to\int_K\psi\,d\mu$ while
$\|d(x,K)^{|\alpha|}\,D_x ^\alpha \phi_\epsilon(x)\|_{L^1}\le  \|
\epsilon^{|\alpha|}\,  D_x^\alpha  \phi_\epsilon(x)\|_{L^1}\to0$
as $\epsilon\to0$ (when $\alpha=0$ one uses the fact that
$|K|=0$). 

Thus, (3.14) and (3.16) show that
 $$ \int_K\psi\,d\mu=0,\quad\psi\in\ccinf(B(0,a)), $$ which implies that the
same conclusion holds for any continuous function $\psi$ on $K$ (first extend
$\psi$ to a compactly supported function on $B(0,a)$ and then approximate the
extension by test functions). Thus the total variation $|\mu|(K)$ of $\mu$ on
$K$ is zero and  by the regularity of $\mu$ it follows that $|\mu|(F')=0$
whenever $F'\subset F$ is a Borel set with  $|F'|=0$. This proves that $\mu_F$
is absolutely continuous with respect to Lebesgue measure.

We now consider the set
$$
F_0=\{x  \in B(0,a):\quad \exists \epsilon>0 : b_j(x,t)=0,\, \forall t\in[0,\epsilon],\, j=0,\dots,n\}
$$
which is a countable union of closed sets
$$
F_\epsilon=\{x  \in B(0,a):\quad  b_j(x,t)=0,\, \forall t\in[0,\epsilon],\, j=0,\dots,n\}
$$
to which we can apply Lemma 3.3 and conclude that
$\mu_{F_\epsilon}$ is absolutely continuous with respect to
Lebesgue measure. Thus, $\mu_{F_0}$ is also absolutely
continuous with respect to Lebesgue measure and the Radon-Nikodym
theorem implies that there exists $g\in L^1_{\text{loc}}(B(0,a))$ such that $$ \mu_{F_0}(X)=\int_X\,
g(x)\,dx,\quad X\subset B(0,a)\text { a Borel set.} $$
\enddemo
The results proved sofar immediately imply Theorem 1.1:

\demo{Proof of Theorem 1.1} We may assume that the vector field
has the form $$ L={\partial\over\partial
t}+ib(x,t){\partial\over\partial x} $$ where $b(x,t)$ is real and
smooth on a neighborhood of $U=B(0,a)\times(-T,T)$ of the origin
in $\erre^{2}$ with $B(0,a)=\{x\in\erre:\,\,\,|x|<a\}$. Since the
trace $bf$ is a measure, by the Radon-Nikodym theorem, we may
write $$ bf=g+\mu $$ where $g$ is a locally integrable function
and $\mu$ is a measure supported on a set $E$ of Lebesgue measure
zero.
Suppose  $x_0$ is a point for which we can find a sequence
$t_j$ converging to $0$ with $b(x_0,t_j)\neq 0$. Let $Z(x,t)$ be a
first integral satisfying $Z(x_0,0)=0$, and $Z_x(x_0,0)=1$. If Im
$ Z_t(x_0,0)\neq 0$, then $L$ will be elliptic in a neighborhood
of $(x_0,0)$ and so by the classical F. and M. Riesz theorem, we
can conclude that $bf$ is absolutely continuous near $(x_0,0)$.
Otherwise, the proof of Theorem 3.1 shows that the FBI transform
with this $Z$ as a first integral and arbitrarily large $\kappa$
decays exponentially in a complex conic neighborhood of
$(x_0,\xi_0)$, for some nonzero covector. By Theorem 2.2 in [BCT],
it follows that near the point $x_0$, modulo a smooth
nonvanishing multiple, the trace $bf$ is the weak boundary value of
a holomorphic function $F$ defined on a side of the curve
$x\longmapsto Z(x,0)$.  But then,  again by the classical F. and M. Riesz
theorem, $bf$ is locally integrable near $x_0$, that is,
$x_0\notin E$. Hence the set $E$ is contained in the  set
$$
F_0=\{x  \in B(0,a):\quad \exists \epsilon>0 : b_j(x,t)=0, \forall t\in[0,\epsilon], j=0,\dots,n\}.
$$
But we already observed that the restriction of $bf$ to $F_0$ is absolutely continuous
with respect to Lebesgue measure which implies that  $\mu$ is  zero.
\enddemo
\proclaim{Remarks}
{\rm
\roster
\item
In the preceding proof, instead of using Theorem
2.2 in [BCT], we can use Lemma 3.3 and Theorem 3.1 in this paper together with
Theorem 1.4 in [B]. However, since this latter theorem in [B]
uses a deep theorem of Uchiyama on the characterization of the
real Hardy space, we chose to present a simpler argument.
\item In the proof of part ii) of Theorem 3.1 we showed how to construct ---under the hypothesis 
that $L$ is proportional to a real vector field along an integral curve of
that vector--- a $C^1$ solution such that its trace has a full wave front set
at the origin. An obvious modification of the proof yields $C^k$ solutions with
the same property for any $k=1,2,\dots$ and the question arises whether it
would be possible to take $k=\infty$. Cleary, this is not true in general
because singularties of the trace may propagate to the interior, as it is easy
to check with the simple example  $L=\partial_t$ where the solutions $Lf=0$ are
funcions of $x$ alone. We will return to this matter in Example 4.3 in the
next section. \endroster} \endproclaim

\heading 4. Examples and applications
\endheading
 We begin here with a lemma which shows that when the vector field is locally
integrable, then a solution is determined by its trace. More precisely, we
have :

\proclaim{Lemma 4.1} Suppose $X$, $U$, $L$ and $f$ are as in Lemma 1.2
and assume in addition that $Lf=0$ in $U_+$, $L$ is locally integrable
in $U$ and that the trace $bf=0$ in $X$. Then $f\equiv 0$ in a neighborhood of $X\times\{0\}$ in $U_+\cup(X\times\{0\})$.
\endproclaim
\demo{Proof} Estimate (1.2) in section 1 allows us to define a
distribution $h$ in $U$ by

$$
\langle h,\psi(x,t)\rangle=\int_0^T\int_Xf(x,t)\psi(x,t)dxdt
$$
 We will show that $Lh=0$ in $U$. Since $h=f$ when $t>0$ and $h=0$ when
$t<0$, we need only show that $h$ is a solution near $t=0$. Suppose then
$\phi(x)$ and $\psi(t)$ are smooth functions of compact support and
$\psi(T)=0$. We have:
$$
\align
\langle Lh,\phi(x)\psi(t)\rangle
    &=\langle h,L^*(\phi(x)\psi(t))\rangle =-\lim_{\epsilon\to
     0^+}\int_{\epsilon}^T\int_Xf(x,t)L^*(\phi(x)\psi(t))dxdt\\
    &=-\lim_{\epsilon \to  0^+} \int_{\epsilon}^T\int_X(Lf(x,t))\phi(x)\psi(t)dxdt\\
    &\quad +\left(\int_Xf(x,T)\phi(x)dx\right)\psi(T)-\lim_{\epsilon \to
    0^+}\int_Xf(x,\epsilon)\phi(x)dx\psi(\epsilon)\\
    &=0.
\endalign
$$
 Note that the second equality above is justified by estimate (1.3).
Thus $Lh=0$ in $U$, and since the trace of $h$ on a noncharacteristic
hypersurface is zero, by a well known theorem of uniqueness for locally
integrable vector fields (see [BT]), it follows that $h\equiv 0$.
\enddemo

We will next apply the method of proof of Theorem 3.1 to present an
example of a nonanalytic tube vector field in the plane which exhibits
an interesting property: the trace of any $C^1$ solution of tempered
growth is real analytic and extends as a smooth solution in a full
neighborhood of points on the boundary.

 \example{Example 4.2} Let
$e(t)=\exp(-1/t^2)$ and set $$
\aligned
\Phi(x_1,x_2,t)&=(\Phi_1,\Phi_2)=
    e(t)\left(\cos(t^{-1}),\sin(t^{-1})\right),\\
            L&={\partial\over\partial t}-
            i\Phi_{1t}{\partial\over\partial x_1}-
            i\Phi_{2t}{\partial\over\partial x_2}\\
            Z&=(x_1+i\Phi_1,x_2+i\Phi_2).\\
\endaligned
$$
 Then, for every unit vector $v$ in $\erre^2$ there is a sequence
satisfying (1), (2) and (3) of Lemma 3.2. Since the origin can be
replaced by any point in the argument, it follows from the proof of
Theorem 3.1 and the fact that $Z_j(x,0)=x_j$ that the weak trace of any
$C^1$ solution of $Lf=0$, $t>0$, with tempered growth as $t\searrow0$,
has to be real analytic. Say $f(x,0)=\sum a_{\alpha}x^{\alpha}$ , where
the power series converges in some neighborhood of the origin. Let
$H(x,t)=\sum a_{\alpha}Z(x,t)^{\alpha}$. Then $ LH=0$ in a neighborhood
of the origin in the plane and by Lemma 4.1, $H$ agrees with $f$ in the
region $t>0$.

\endexample

\example{Example 4.3}
Consider  the vector field  in the plane given by
$$
 L=\frac{\partial}{\partial
 t}-\frac{2ixt}{1+it^2}\frac{\partial}{\partial x}\tag 4.1
$$
The $t$-axis is an elliptic submainifold for $L$ (see [HT] for the definition) 
and by the main result proved in [HT] this axis propagates analyticity for
solutions of the homogeneous equation. That is, if $Lu=0$ on a neighborhood
$\Omega$ of the origin and $u$ is analytic at some point $(0,t_0) \in\Omega$,
then $u$ is analytic at every point $(0,t)\in\Omega$. However,  Treves showed
[T] that the $t$-axis does not propagate smoothnes: 
 there is a solution $u$
of $L$ which is continuous in the plane, smooth off  $\{t=0\}$, but the trace
$u(x,0)$ is not
 smooth at the origin. The existence of such a solution was
proved in a
 nonconstructive fashion using a Baire category argument.
 Here we
wish to construct a solution $h$ with the additional property that the
$C^{\infty}$ wave front set of the trace $h(x,0)$ at the origin contains both
directions $1$ and $-1$.
\newline
 Observe that the function $Z(x,t)=x(1+it^2)$
is a first integral of $L$ in the
 plane. For each $k$ a positive integer,
define
 $$
  W_k(x,t)=\left (Z(x,t)^2-\frac{1}{k^2}\right )^{\frac{1}{2}}
 $$
where we take the principal branch of the square root off the negative $y$
axis. Note that $W_k$ is continuous in the plane and smooth in $\Bbb C $
except at $(\frac{1}{k},0)$ and $(-\frac{1}{k},0)$. Let
 $$
  h(x,t)=\sum_{k=1}^{\infty}\frac{1}{3^k}W_k(x,t)
$$
The function $h$ is clearly continuous everywhere since the series converges
absolutely on compact sets. We will show that $h$ is smooth when $t\neq 0$.
Let
$$
 g_k(x,t)=Z(x,t)^2-\frac{1}{k^2}
$$
Fix two positive numbers $M>\delta$ and consider the size of $g_k$ in the
region $\delta \leq |t| \leq M$. We have
$$
 |g_k(x,t)|^2=(1+2t^4+t^8)x^4+(\frac{2t^4-2}{k^2})x^2+\frac{1}{k^4}
$$
It follows that there exists a constant $C(M,\delta) >0$ such that
$$
 |g_k(x,t)|\geq \frac{C(M,\delta)}{k^2} \tag 4.2
$$
whenever $\delta \leq |t| \leq M$. We will now show that when $t\not=0$ and
for $n$ a positive integer,
$$
D^nW_k(x,t)=\sum_{j=0}^{\left[\frac{n}{2}\right]}A_jg_k^{\frac{1}{2}-n+j}(Dg_k)^{n-2j}(D^2g_k)^j
$$
where $D=\frac{\partial}{\partial x}$, $[y]$ denotes the greatest integer less
than or equal to $y$, and the $A_j$ are constants depending only on $n$. In
particular, these constants do not depend on $k$. We will prove this assertion
by inducting on $n$. When $n=1$, the formula is valid since
$$
DW_k(x,t)=\frac{1}{2}(Dg_k)g_k^{-\frac{1}{2}}
$$
Assume the assertion for some $n$ and apply $D$ to both sides of the formula.
When $n$ is odd, since $D^3g_k=0$, we get
$$
\aligned
D^{n+1}W_k&=\sum_{j=0}^{\left[\frac{n}{2}\right]}(\frac{1}{2}-n+j)A_jg_k^{\frac{1}{2}-n-1+j}(Dg_k)^{n+1-2j}(D^2g_k)^j  \\
&+\sum_{j=0}^{[\frac{n}{2}]}A_jg_k^{\frac{1}{2}-n+j}(Dg_k)^{n-1-2j}(D^2g_k)^{j+1}
\endaligned
 $$
In the second sum replace $j$ by $j+1$ and observe that since
$n$ is odd,
$$
\left[\frac{n}{2}\right]+1=\left[\frac{n+1}{2}\right]
$$
We are then led to
$$
D^{n+1}W_k(x,t)=\sum_{j=0}^{\left[\frac{n+1}{2}\right]}B_jg_k^{\frac{1}{2}-n-1+j}(Dg_k)^{n+1-2j}(D^2g_k)^j
$$
 for some $B_j$ depending only on the $A_l$ and $n$, and hence
independent of $k$. When $n$ is even, observe that in the second sum of
$D^{n+1}$, the index $j$ goes only upto $[n/2]-1$. Hence in
this second sum if we replace $j$ by $j+1$ and observe that
$$
\left[\frac{n}{2}\right]=\left[\frac{n+1}{2}\right],
$$
we get an expression for $D^{n+1}W_k$
as required.\newline The expression for $D^nW_k(x,t)$ together with the
lower bound $(4.2)$ show that the series $$
\sum_{k=1}^{\infty}\frac{1}{3^k}D_x^nW_k(x,t) $$ converges absolutely
for $x$ in a compact set and $0<\delta\leq |t|\leq M$.  Thus
$D_x^nh(x,t)$ exists for all $n$ when $t\neq 0$. Next note that since
$LW_k(x,t)=0$ when $t\neq 0$, we have $$
D_x^nD_tW_k=\sum_{j=0}^{n+1}P_j(x,t)D_x^jW_k, $$ for some smooth $P_j$.
Therefore using what was already proved, we see that for any $n$,
$D_x^nD_th(x,t)$ exists and $Lh=0$ when $t\neq 0$. We can now iterate
by differentiating the equation $Lh=0$ to conclude that $h$ is smooth
when $t\neq 0$. Finally, note that the trace $$
h(x,0)=\frac{1}{3^k}W_k(x,0)+E_k(x) $$ where $E_k$ is $C^1$ at the
points $|x|=\frac{1}{k}$. Hence $h(x,0)$ is not $C^1$ in any
neighborhood of the origin. Moreover, since $h(x,0)=h(-x,0)$, both $1$
and $-1$ are in the wave front set of $h(x,0)$ at the origin.
\endexample
\proclaim{Remarks}
{\rm
\roster
\item
 In the preceding example, if we restrict the solution $h(x,t)$ to the domain $D=\{(x,t):\,\,\, t>x^2\}$ bounded by the parabola $t=x^2$ we see that its boundary value $bh(x)=h(x,x^2)$  is smooth except at the origin where its wave front set contains both directions.
\item
The geometric background behind Example 4.3 is as follows: the vector field
$L$ is a Mizohata vector field for $x>0$ and a conjugate Mizohata vector field
for $x<0$. The wave front set of the trace of a smooth solution of the
homogeneous Mizohata equation defined on the upper plane necessarily lacks the
direction $\xi<0$ while in the case of the conjugate Mizohata vector field the
missing direction is $\xi>0$. However, approaching the origin from both sides
a singularity that contains the two microlocal directions can be produced.
\roster

}
\endproclaim


\Refs
\widestnumber\key{BCT} \refstyle{A}
\ref
\key BT
\by M. S.
Baouendi and F. Treves
\paper A property of the functions and
distributions annihilated by a locally integrable system of
complex vector fields
\jour Ann. of Math.
\vol 113
\yr 1981
\pages
387--421
\endref
\ref
\key BCT
\by M. S. Baouendi, C. H. Chang,
and F. Treves
\paper Microlocal hypo-analyticity and extension of
CR functions
\jour J. Diff. Geom.
\vol 18
\yr 1983
\pages 331--391
\endref

\ref
\key BC
\by S. Berhanu and S. Chanillo
\paper Boundedness of the FBI transform on Sobolev spaces and propagation of
singularities
\jour Commun. in PDE
\vol 16 (10)
\yr 1991
\pages  1665--1686
\endref

\ref
\key B
\by R. G. M. Brummelhuis
\paper A Microlocal F.and M.Riesz Theorem with applications
\jour Revisita Matematica Iberoamericana
\vol 5
\yr 1989
\pages 21--36
\endref

\ref 
\key  HT
\by  N. Hanges and F. Treves
\paper Propagation of holomorphic extendability of CR functions
\jour Math. Ann.
\vol 263
\yr 1983 
\pages 157--1 77
\endref 

\ref
\key H
\by J. Hounie
 \paper Globally hypoelliptic vector fields on compact surfaces
\jour Comm. Partial Differential Equations
\vol 7
\yr 1982
\pages 343--370
\endref

\ref
\key NT
\by L. Nirenberg and F. Treves
\paper Solvability of a first order linear partial differential equation
\jour Comm. Pure Appl. Math.
\vol 16
\yr 1963
\pages 331--351
\endref

\ref
\key S
\by J. Sjostrand
\paper Singularit\'es analytiques
microlocales
\jour Ast\'erisque
\vol 95
\yr 1982
\endref

\ref
\key T
\by F. Treves
\book Hypo-analytic structures, local theory
\publ Princeton University Press
\yr 1992
\endref

\ref 
\key  U
\by A. Uchiyama 
\paper A constructive proof of the Fefferman-Stein decomposition of
$BMO(\erre^n)$
\jour Acta Math.
\vol 148 
\yr 1982 
\pages 215--241 
\endref 

\endRefs
\enddocument